\documentclass[10pt]{amsart}

\newenvironment{prop}{\medskip\noindent {\bf Proposition.}\it}{\rm\par\medskip}
\newenvironment{myproof}{\medskip\noindent{\it Proof\/}.}{$\;\square$
                         \par\medskip}
\newenvironment{mydef}{\medskip\noindent {\bf Definition.}\it}{\rm\par\medskip}
\newenvironment{remark}{\medskip\noindent {\bf Remark.}}{\par\medskip}

\newcommand{\cp}[3]{Comp. Math. #1 (#2), #3}
\newcommand{\dmj}[3]{Duke Math. J. #1 (#2), #3}
\newcommand{\faa}[3]{Funct. Anal. Appl. #1 (#2), #3}
\newcommand{\imrn}[3]{Int. Math. Res. Notices #1 (#2), #3}
\newcommand{\jpamg}[3]{J. Phys. A: Math. Gen. #1 (#2), #3}
\newcommand{\np}[3]{Nucl. Phys. #1 (#2), #3}

\newcommand{\tmp}[3]{Theor. Math. Phys. #1 (#2), #3}

\newcommand{\esp}[1]{Elsevier Science Publishers (#1)}

\font\fonta=cmtt10
\def\arob{\hbox{\fonta\char'100}}

\font\fontb=msam10
\def\square{\hbox{\fontb\char'003}}

\newcommand{\Z}{{\mathbb Z}}
\newcommand{\N}{{\mathbb N}}

\newcommand{\C}{{\mathbb C}}

\newcommand{\g}{{\mathfrak{g}}}
\newcommand{\h}{{\mathfrak{h}}}
\newcommand{\hd}{{\mathfrak{h}}^\ast}
\newcommand{\Egln}{E_{\tau,\gamma}(\mathrm{gl}_N)}
\newcommand{\Egltwo}{E_{\tau,\gamma}(\mathrm{gl}_2)}
\newcommand{\gln}{\mathrm{gl}_N}
\newcommand{\End}{{\mathrm{End}}}
\newcommand{\sumi}[1]{\sum_{i=#1}^N}
\newcommand{\sumij}[1]{\sum_{_{\ i \neq j}^{i,j=#1}}^N}
\newcommand{\prodij}{\prod_{j;j\neq i}}
\newcommand{\tr}{\mathrm{tr}}
\newcommand{\LL}{\mathcal{L}}
\newcommand{\T}{\mathcal{T}}
\newcommand{\Rred}{\check{R}}
\newcommand{\Pred}{\check{P}}
\newcommand{\lared}{\check{\lambda}}
\newcommand{\hred}{\check{h}}
\newcommand{\Tred}{\check{\T}}
\newcommand{\psired}{\check{\psi}}
\newcommand{\Gammared}{\check{\Gamma}}
\newcommand{\Lred}{\check{\LL}}
\newcommand{\epsred}{\check{\varepsilon}}
\newcommand{\Bred}{\check{B}}
\newcommand{\vred}{\check{v}}
\newcommand{\gred}{\check{g}}
\newcommand{\Dred}{\check{D}}
\newcommand{\Phired}{\check{\Phi}}
\newcommand{\Wred}{\check{W}}
\newcommand{\Rcalred}{\check{\mathcal{R}}}
\newcommand{\omegared}{\check{\omega}}
\newcommand{\Gam}[1]{\Gammared^{(#1)}}
\newcommand{\Gaminv}[1]{\left(\Gam{#1}\right)^{-1}}
\newcommand{\RR}[2]{\Rred^{(#1 #2)}}

\begin{document}

\title[Bethe ansatz for the Ruijsenaars model]
{Algebraic nested Bethe ansatz for the \\ elliptic Ruijsenaars model}
\author{E. Billey}
\maketitle

\begin{center}
{\it D-MATH, ETH-Zentrum, 8092 Z\"urich, Switzerland} \\
{\tt billey$\arob$math.ethz.ch}
\end{center}

\bigskip

\begin{abstract}
The eigenvalues of the elliptic $N$-body Ruijsenaars operator 
are obtained by a dynamical version of the algebraic nested Bethe 
ansatz method. The result is derived by using the 
construction given in~\cite{FeVa97}, 
where the Ruijsenaars operator was obtained as the
transfer matrix associated to the symmetric power of the vector 
representation of the elliptic quantum group $\Egln$.
\end{abstract}

\bigskip

\section{Introduction}

The goal of this paper is to solve the elliptic $N$-body Ruijsenaars 
model by the Bethe ansatz method. It was done in~\cite{FeVa96} for 
$N=2$, as a special case of the diagonalization 
of the transfer matrix of particular highest 
weight modules associated to the elliptic quantum group $\Egltwo$.
The result was achieved by a dynamical generalization of the algebraic Bethe 
ansatz (see~\cite{SkTaFa79, Fa82} for the description of 
the algebraic formulation of the Bethe ansatz). For general $N$, the 
Ruijsenaars operator can be constructed as the transfer matrix $\T(z)$ 
associated to the symmetric power of the vector representation of 
$\Egln$~\cite{FeVa97}. 
Since $\T(z)$ is a trace over an $N$-dimensional space, we diagonalize it 
by the nested version of the Bethe ansatz~\cite{KuRe83}. 

The Ruijsenaars operator is a difference operator which is a 
q-deformation of the Calogero differential operator. In~\cite{FeVa95a} 
the solution of the elliptic $N$-body Calogero model (in the Bethe 
ansatz form) was obtained from the integral representation of solutions 
of the elliptic Knizhnik--Zamolodchikov--Bernard equations by applying 
the stationary phase method~\cite{EtKi93a, EtKi94, ReVa94}. 
Similarly, there is a link between the Ruijsenaars model and the 
qKZB difference equations. The eigenfunctions of the trigonometric 
Ruijsenaars operator are given in~\cite{FeVa95b} in the form of 
Etingof--Kirillov traces of intertwining operators~\cite{EtKi93b}. As for 
the integral representation of solutions of the qKZB equations, 
it is known for $\mathrm{sl}_2$ only~\cite{FetaVa96}. 

The structure of the paper is as follows. 
In section~\ref{2}, we review  
the construction of commuting transfer matrices associated to the 
representations of the elliptic quantum group $\Egln$. 
In section~\ref{3}, we recall how the Ruijsenaars operator 
is related to the transfer matrix $\T(z)$ associated to the symmetric 
power of the vector representation of $\Egln$. 
In section~\ref{4}, we explain the idea of the dynamical version of 
the algebraic nested Bethe ansatz. Finally in sections~\ref{5} and~\ref{6}, 
we apply explicitly the Bethe ansatz method; we write down the Bethe 
equations and the eigenvalues of the transfer matrix $\T(z)$. 

\section{The elliptic quantum group associated to $\gln$}

\label{2}

We recall the definition of the elliptic quantum group $\Egln$, 
following~\cite{FeVa97}. 
Let $\g$ be a Lie algebra, $\h$ a Cartan subalgebra of $\g$ and $\hd$ 
its dual space. Let $W$ be a finite dimensional diagonalizable 
$\h$-module, {\it i.e.} a complex finite dimensional vector space 
with a weight decomposition 
$W=\oplus_{\mu \in \hd} W[\mu]$ such that $\h$ acts on $W[\mu]$ by 
$x w = \mu(x) w$ for $x \in \h, w \in W[\mu]$.

The starting point to define representations of elliptic quantum 
groups is 
an $R$-matrix $R(z,\lambda) \in \End(W \otimes W)$ depending on two 
parameters $z \in \C, \lambda \in \hd$ and solution of the dynamical 
Yang--Baxter equation 
\begin{eqnarray}
\label{DYB} 
\lefteqn{R^{(12)}(z_1-z_2, \lambda - \gamma h^{(3)}) \ 
R^{(13)}(z_1-z_3, \lambda) \ 
R^{(23)}(z_2-z_3, \lambda - \gamma h^{(1)})} \\ 
\nonumber
& & = \ R^{(23)}(z_2-z_3, \lambda) \ 
R^{(13)}(z_1-z_3, \lambda - \gamma h^{(2)}) \ 
R^{(12)}(z_1-z_2, \lambda).  
\end{eqnarray}
In this equation $\gamma$ is a generic complex parameter, and we 
use the standard notation
\begin{equation}
\label{shift}
f(\lambda-\gamma h) \ w \ = \ f(\lambda-\gamma \mu) \ w \ \ \ \ \mathrm{if} 
\ w \in W[\mu]
\end{equation}
for any complex-valued function $f$ of $\lambda$. 

In what follows we are interested in the case $\g=\gln$, and $\h$ the 
algebra of diagonal complex $N \times N$ matrix, acting on the vector 
representation $V=\C^N$ of $\gln$ with standard basis 
$(e_j)_{j=1,\ldots,N}$. We identify $\h$ and $\hd$ via the trace, and 
with $\C^N$ via the orthonormal basis $(\omega_j=e_{jj})_{j=1,\ldots,N}$, 
denoting  $e_{jk}$ the $N \times N$ matrix acting on the 
standard basis by $e_{jk} e_l = \delta_{kl} e_j$. The weight spaces of 
the $\h$-module $V$ are $V[\omega_j]=\C e_j$. In this case the 
notation~(\ref{shift}) reads explicitly
\begin{displaymath}
f(\lambda-\gamma h) \ e_k \ = \ \Gamma_k f(\lambda) \ e_k
\end{displaymath}
where $\Gamma_k$ is the shift operator by $-\gamma \omega_k$: if $f$ 
is a complex-valued function of $\lambda$, 
\begin{displaymath}
\Gamma_k f(\lambda) \ = \ f(\lambda-\gamma\omega_k) \ = \ 
f(\lambda_1, \ldots, \lambda_{k-1}, 
\lambda_k-\gamma, \lambda_{k+1}, \ldots, \lambda_N).
\end{displaymath}

Let $R(z,\lambda) \in \End(\C^N \otimes \C^N)$ be the elliptic solution of 
the dynamical Yang--Baxter equation~(\ref{DYB}) given by the formula
\begin{equation}
\label{Rmatrix}
R(z,\lambda) \ = \ \sumi{1} e_{ii} \otimes e_{ii} + \sumij{1} 
\alpha(z,\lambda_i - \lambda_j) \ e_{ii} \otimes e_{jj} 
+ \sumij{1} \beta(z,\lambda_i - \lambda_j) \ e_{ij} \otimes e_{ji}, 
\end{equation}
where 
\begin{eqnarray*}
\alpha(z,\lambda) & = &
\frac{\theta(z) \theta(\lambda+\gamma)}{\theta(z-\gamma) \theta(\lambda)}, 
\ \ \ \ \beta(z,\lambda) = 
- \frac{\theta(z+\lambda) \theta(\gamma)}{\theta(z-\gamma) \theta(\lambda)}. 
\end{eqnarray*}
$\theta$ is the Jacobi's first theta function:
\begin{displaymath}
\theta(z) = - \sum_{j \in \Z} \exp \left[ \pi i (j+ \frac{1}{2})^2 \tau 
+ 2 \pi i (j+ \frac{1}{2}) (z+\frac{1}{2}) \right],
\end{displaymath}
where $\tau$ is a complex parameter with $\mathrm{Im}(\tau) > 0$. 
$\theta$ is analytic with zeroes $(n+\tau m), n,m \in \Z$. It 
satisfies
\begin{displaymath}
\theta(-z) \ = \ - \ \theta(z) \ = \ \theta(z+1) \ \ \ \ \mathrm{and} \ \ \ \ 
\theta(z+\tau) \ = \ - \ \theta(z) \ \mathrm{e}^{-i \pi (2z+\tau)}.
\end{displaymath}

The $R$-matrix~(\ref{Rmatrix}) has the following properties:
\begin{eqnarray}
& & R^{(12)}(z,\lambda) \ R^{(21)}(-z,\lambda) \ = \ \mathrm{Id}, \\
& & R^{(12)}(0,\lambda) \ = \ P^{(12)} \ = \ \sum_{i,j=1}^N 
e_{ij} \otimes e_{ji}, \\
\label{commRh}
& & \relax [ R^{(12)}(z,\lambda) , x \otimes \mathrm{Id} + 
\mathrm{Id} \otimes x ] \ = \ 0, \ \ \ \ \forall x \in \h, \\
\label{shiftR}
& & \relax [ \Gamma^{(1)} \Gamma^{(2)}, R^{(12)}(z,\lambda) ] \ = \ 0, 
\end{eqnarray}
where $\Gamma$ is the diagonal $N \times N$ matrix $\mathrm{diag}
(\Gamma_i)_{i=1,\ldots,N}$. 

A representation of the elliptic quantum group $E=\Egln$ is by definition
a pair $(W,L)$ where $W$ is a finite-dimensional diagonalizable 
$\h$-module and $L(z,\lambda)$
is a meromorphic function with values in $\End_\h(\C^N\otimes W)$,
obeying the relation 
\begin{eqnarray*}
\lefteqn{R^{(12)}(z_1-z_2,\lambda-\gamma h^{(3)}) \ L^{(13)}(z_1,\lambda) \ 
L^{(23)}(z_2,\lambda-\gamma h^{(1)})} \\
& & = \ L^{(23)}(z_2,\lambda) \ 
L^{(13)}(z_1,\lambda-\gamma h^{(2)}) \ R^{(12)}(z_1-z_2,\lambda)
\end{eqnarray*}
and commuting with the action of $\h$:
\begin{displaymath}
[ L(z,\lambda) , x \otimes \mathrm{Id} + \mathrm{Id} \otimes x ] \ = \ 0, 
\ \ \ \ \forall x \in \h.
\end{displaymath}
An $E$-submodule of an  $E$-module $(W,L)$ is a pair $(W',L')$
where $W'$ is an $\h$-submodule of $W$ such that $\C^N\otimes W'$
is invariant under the action of $L(z,\lambda)$, and
$L'$ is the restriction of $L$ to this invariant
subspace. $E$-submodules are $E$-modules.

The basic example of an $E$-module is
$(\C^N,L)$ with $L(z,\lambda)=R(z-w,\lambda)$. It is called the 
vector representation with evaluation point $w$ and is denoted 
by $V(w)$.

Other modules can be obtained by taking tensor
products: if $(W_1,L_1)$ and $(W_2,L_2)$ are $E$-modules,
then also $(W_1\otimes W_2,L)$, with an $\h$-module structure 
$x (w_1 \otimes w_2) = x w_1 \otimes w_2 + w_1 \otimes x w_2$ 
and a $L$-operator $L(z,\lambda)=
L_1^{(12)}(z,\lambda-\gamma h^{(3)}) \ L_2^{(13)}(z,\lambda)$.

It is useful for what follows to introduce what we call the 
Lax operator associated to the $E$-module $(W,L)$. It is 
an $N \times N$ matrix, with elements which are 
operators acting on the space of meromorphic functions of $\lambda\in\h^*$ 
with values in $W$. It is defined by the formula 
\begin{equation}
\label{defL}
\LL(z) \ w(\lambda) \ = \ L^{(12)}(z,\lambda) \ \Gamma^{(1)} \ 
w^{(2)}(\lambda). 
\end{equation} 
More explicitly, let us introduce matrix elements
by $L(z,\lambda) \ e_j \otimes w = \sum_i e_i \otimes L_{ij}(z,\lambda) \ w$ 
(and the same for $\LL$). The elements $\LL_{ij}(z)$ of $\LL(z)$ 
act the following way:
\begin{displaymath}
\LL_{ij}(z) \ w(\lambda) \ = \ L_{ij}(z,\lambda) \ w(\lambda-\gamma\omega_j).
\end{displaymath}

Using the property~(\ref{shiftR}) of the $R$-matrix, we see easily that 
$\LL$ satisfies the commutation relation 
\begin{equation}
\label{RLL}
R^{(12)}(z_1-z_2, \lambda-\gamma h^{(3)}) \ \LL^{(13)}(z_1) \ \LL^{(23)}(z_2)
= \LL^{(23)}(z_2) \ \LL^{(13)}(z_1) \ 
R^{(12)}(z_1-z_2, \lambda).
\end{equation}

The trace of the Lax operator 
$\tr^{(1)} \LL^{(13)}(z) = \sum_i \LL_{ii}(z)$ leaves the zero weight 
subspace $W[0]$ of $W$ invariant. The transfer matrix associated to 
the $E$-module $(W,L)$ is, by definition,
\begin{displaymath}
\T(z) = \left[ \tr^{(1)} \LL^{(13)}(z) \right]_{W[0]}.
\end{displaymath}
As a consequence of relation~(\ref{RLL}), the transfer matrices 
commute for different values of the spectral parameters. 

\section{The Ruijsenaars operator}

\label{3}

The elliptic Ruijsenaars operator is (up to a conjugation by a function) a  
difference operator acting on functions of $\lambda \in \hd$:
\begin{displaymath}
M \ = \ \sumi{1} \prodij 
\frac{\theta(\lambda_i-\lambda_j+\ell\gamma)}{\theta(\lambda_i-\lambda_j)} 
\ \Gamma_i, \ \ \ \ \mathrm{with \ a \ coupling \ constant} \ \ell \in \N.
\end{displaymath}

It can be obtained as the transfer matrix 
associated to a particular $E$-module  that we now introduce. 

Let $n \in \N$. Let $V^{\otimes n}(0)$ denote the $E$-module 
$V(0) \otimes V(\gamma) \otimes \cdots \otimes V(\gamma(n-1))$. 
$V^{\otimes n}(0)$ is the pair $(W,L)$ with $W=(\C^N)^{\otimes n}= 
\C^N \otimes \cdots \otimes \C^N$, and $L$ is the operator on 
$\C^N \otimes (\C^N)^{\otimes n}$ given by the formula
\begin{displaymath}
L(z,\lambda) \ = \ R^{(01)}\Big(z,\lambda-\gamma\sum_{j=2}^n h^{(j)}\Big) \ 
R^{(02)}\Big(z-\gamma,\lambda-\gamma\sum_{j=3}^n h^{(j)}\Big) \ \cdots \ 
R^{(0n)}\Big(z-\gamma(n-1),\lambda\Big).
\end{displaymath}

We denote by $S^n(\C^N)$ the space of symmetric tensors 
of $(\C^N)^{\otimes n}$, {\it i.e.} the 
subspace of $(\C^N)^{\otimes n}$ invariant under the action of the 
symmetric group $S_n$. 
It is proved in~\cite{FeVa97} that $S^n(\C^N)$ is an $E$-submodule of 
$V^{\otimes n}(0)$. It is called the $n$th symmetric power of 
the vector representation (with evaluation point $0$) and is denoted by 
$S^nV(0)$. The zero weight subspace of this module is trivial unless $n$ 
is a multiple of $N$; if we take $n=N\ell$, the zero weight 
subspace is one dimensional and is spanned by the sum of the tensors 
$e_{i_1} \otimes \cdots \otimes e_{i_{N\ell}}$ over all sequences $(i_j)$ 
such that each integer between 1 and N occurs precisely $\ell$ 
times\footnote{By extension we call zero weight a weight which 
is a multiple of $\omega=\sum_{i=1}^N \omega_i$. Since the $R$-matrix 
and the Lax operator depend only on differences $(\lambda_i-\lambda_j)$, 
adding the same constant to each $\lambda_i$ does not change our results.}. 
Let $\T(z)$ be the transfer matrix associated to $S^{N\ell}V(0)$. If 
we identify the zero weight subspace of $S^{N\ell}V(0)$ with $\C$, we have 
\begin{equation}
\label{Ruij}
\T(z) \ = \ \frac{\theta(z-\gamma\ell)}{\theta(z-\gamma N\ell)} M. 
\end{equation}
Thanks to this result, we can solve the Ruijsenaars model by 
applying the Bethe ansatz method to the transfer matrix $\T(z)$. 

Let us remark (although we do not use this result in what follows) 
that $M$ belongs to a family of $N$ commuting difference operators 
\begin{displaymath}
M_n \ = \ \sum_{I;|I|=n} \ \prod_{^{i \in I}_{j \not\in I}} 
\frac{\theta(\lambda_i-\lambda_j+\ell\gamma)}{\theta(\lambda_i-\lambda_j)} 
\ \prod_{i \in I} \Gamma_i, \ \ \ \ \ \ n=1,\ldots,N
\end{displaymath}
(where $|I|$ denotes the cardinality of a subset $I$ of $\{1,\ldots,N\}$). 
Each operator $M_n$ can be constructed by considering transfer matrices 
associated to $E$-modules obtained as symmetric and exterior powers of 
the vector representation of $\Egln$. 

\section{Setting up of the Bethe ansatz}

\label{4}

Let us recall that the algebraic Bethe ansatz is a method for 
diagonalizing transfer matrices obtained as a trace of a $2 \times 2$ 
Lax matrix with a spectral parameter. The nested Bethe ansatz is a 
generalization of the Bethe ansatz to the case of Lax matrices with 
an $N$-dimensional auxiliary space, for $N>2$; it is achieved in 
(N-1) steps, each step resulting in reducing by 1 the dimension of 
the auxiliary space. 

Here we apply a dynamical version of the nested Bethe ansatz: the 
Lax operator is a matrix on an auxiliary $N$ dimensional space which 
elements are difference operators acting on functions of $N$  
parameters $(\lambda_1,\ldots,\lambda_N)$ with values in some 
vector space.

The idea consists in looking for eigenstates $\psi$ obtained 
by applying ``creation operators'' $B(t_k)$ to a reference state $v$: 
$ \psi(t_1, \ldots, t_m) = B(t_1) \cdots B(t_m) v$. The $B(t_k)$'s  
are particular elements of the Lax matrix $\LL$, and the 
reference state $v$ is chosen to be an obvious eigenstate of the 
elements of the Lax matrix. The problem amounts to 
finding the conditions on the parameters $t_k$. More precisely, 
in the case $N=2$, $B(t)$ is the operator $\LL_{12}(t)$. In the case $N>2$, 
it is a little more complicated: the operators $B(t)$ belong to 
the set $(B_j(t)=\LL_{1j})_{j=2, \ldots,N}$, so that 
$\psi$ is roughly of the form $B_{j_1}(t_1) \cdots B_{j_m}(t_m) v$ 
and we have to choose not only the right parameters $t_k$, but also 
the indices $j_k$. 

Let us now go into details. We want to diagonalize the transfer matrix 
\begin{equation}
\label{genT} 
\T(z) \ = \ \tr^{(0)} \LL(z) 
\end{equation}
where 
\begin{equation}
\label{genL}
\LL(z) \ = \ 
R^{(01)}\Big(z,\lambda-\gamma\sum_{j=2}^{N\ell} h^{(j)}\Big) 
\ \cdots \ 
R^{(0,N\ell)}\Big(z-\gamma(N\ell-1),\lambda\Big) \ \Gamma^{(0)}. 
\end{equation} 
This Lax operator is a matrix on an auxiliary $N$ dimensional 
space (denoted by the index $(0)$) with elements which are difference 
operators acting on functions of $\lambda=(\lambda_1,\ldots,\lambda_N)$ 
with values in $W=(\C^N)^{\otimes N\ell}$. We denote these elements by:
\begin{displaymath}
\LL(z) \ = \ \left( 
\begin{array}{cc}
A(z) & B_j(z) \\ C_i(z) & D_{ij}(z)
\end{array}
\right)_{i,j=2, \ldots,N}.
\end{displaymath}
We denote by $B$ the lign vector $(B_2, \ldots, B_N)$ 
and by $D$ the matrix $(D_{ij})_{i,j=2, \ldots,N}$. 

\subsection{The form of the Bethe ansatz}

Now let us explain more precisely in which form we are looking for the 
eigenstates of $\T(z)$.

First we choose a reference state $v$. We take a joint 
eigenstate of the operators $A$ and $D_{ii}$ which is furthermore 
annihilated by all the $D_{ij}$ for $i \neq j$:
\begin{displaymath}
v \ = \ \underbrace{e_1 \otimes \cdots \otimes e_1}_{N \ell \ 
\mathrm{factors}} \ \in W.
\end{displaymath}

\begin{prop}
The action of the operators $A$ and $D_{ij}$ on the reference state 
is given by the following formulae:
\begin{eqnarray}
\label{Aaction}
& & A(z) \ [g(\lambda) \ v] \ = \ g(\lambda-\gamma\omega_1) \ v, \\
\label{Daction}
& & D(z) \ [g(\lambda) \ v] \ = \ \sum_{i,j=2}^N e_{ij} \ D_{ij}(z)  
\ [ g(\lambda) \ v] \ = \ 
\sumi{2} e_{ii} \ \Phi_i(z,\lambda) \ g(\lambda-\gamma\omega_i)\ v, 
\end{eqnarray}
where $ \displaystyle 
\Phi_i(z,\lambda) \ = \ \frac{\theta(z)}{\theta(z-N\ell\gamma)} 
\frac{\theta(\lambda_i-\lambda_1+N\ell\gamma)}{\theta(\lambda_i-\lambda_1)}$.
\end{prop}

\begin{myproof} 
We can write $R(z,\lambda)=\sum_{i,j=1}^N e_{ij} \otimes 
R_{ij}(z,\lambda)$ with 
\begin{eqnarray*}
& & R_{ii}(z,\lambda) = e_{ii} + \sum_{j;j\neq i} 
\alpha(z,\lambda_i-\lambda_j) \ e_{jj} \ \ \ \ \mathrm{for} \ i=1,\ldots,N, \\
& & R_{ij}(z,\lambda) =  \beta(z,\lambda_i-\lambda_j) \ e_{ji} \ \ \ \ 
\mathrm{for} \ i,j=1,\ldots,N, \ i\neq j.
\end{eqnarray*}
Then we have 
\begin{eqnarray*} 
& & L_{ij}(z,\lambda) \ = \ \sum_{j_1,\ldots,j_{N\ell-1}}^N \ 
R^{(1)}_{i j_1}\Big(z,\lambda-\gamma\sum_{k=2}^{N\ell} h^{(k)}\Big) \ 
R^{(2)}_{j_1 j_2}\Big(z-\gamma,\lambda-\gamma\sum_{k=3}^{N\ell} h^{(k)}\Big) 
\ \cdots \\
& & \hspace*{105pt} \cdots \ R^{(N\ell)}_{j_{N\ell-1} j}
\Big(z-\gamma(N\ell-1),\lambda\Big).
\end{eqnarray*} 

Let us consider $L_{ij}$, for $i,j \geq 2$, $i\neq j$. 
In each term of the sum we have at least one of the 
factors which is of the form $R_{ik}$ for some $k \neq i$. Since 
$R_{ik}$ is proportional to $e_{ki}$, $R_{ik} e_1=0$. Thus 
$D_{ij} \ v = 0$ for $i,j \geq 2$, $i\neq j$. 

Now if we apply $L_{ii}$ to $v$, the only term 
which gives a non zero contribution is the one corresponding to 
$j_1=\cdots=j_{N\ell-1}=i$. Therefore 
\begin{displaymath}
L_{11}(z,\lambda) \ v \ = \ v \ \ \ \ \mathrm{and} \ \ \ \ 
L_{ii}(z,\lambda) \ v \  = \  \Phi_i(z,\lambda) \ v  
\ \ \ \ \mathrm{if} \ i\neq 1, 
\end{displaymath}
with $ \displaystyle \Phi_i(z,\lambda) \ = \ 
\prod_{k=1}^{N\ell} \alpha(z-\gamma(k-1),\lambda_i-\lambda_1+\gamma(N\ell-k))$.
\end{myproof}

Then we look for eigenstates in the form 
\begin{displaymath}
\psi(t_1,\ldots,t_m,\lambda) = \sum_{j_1,\ldots,j_m=2}^N 
B_{j_1}(t_1) \cdots B_{j_m}(t_m) \ v \ g_{j_1,\ldots,j_m}(\lambda),
\end{displaymath}
with some coefficients $g_{j_1,\ldots,j_m}$ which are functions of 
$\lambda$ (and depend implicitely on the parameters $t_k$). 
We impose conditions on these $g_{j_1,\ldots,j_m}$'s 
so as to consider only the states $\psi$ of weight zero. 
$v$ has a weight $[(N-1)\ell,-\ell,\ldots,-\ell]$ in the 
basis $(\omega'_j=e_{jj}-\omega/N)$. Since applying 
the operator $B_j$ to a vector of weight $\mu$ gives a vector of weight 
$\mu-\omega_1+\omega_j$, we have to apply to $v$ 
exactly $\ell$ times each of the $B_j$ in order to get a zero weight vector. 
This means that we must take a state of the form
\begin{equation}
\label{psi}
\psi(t_1,\ldots,t_{(N-1)\ell},\lambda) \ = \ B^{(1)}(t_1) \cdots 
B^{((N-1)\ell)}(t_{(N-1)\ell}) \ v \ g^{(1,\ldots,(N-1)\ell)}(\lambda)
\end{equation}
where $g^{(1,\ldots,(N-1)\ell)}$ is a function of $\lambda$ with values 
in the zero weight subspace of $\Wred=(\C^{N-1})^{\otimes (N-1)\ell}$.

\subsection{The commutation relations}

To evaluate the action of the transfer matrix
\begin{displaymath}
\T(z) \ = \ A(z) + \tr D(z)
\end{displaymath}
on the vector $\psi$, we  have to push $A$ and $D$ to the right 
of the $B$ operators, and we need to write the commutation 
relations~(\ref{RLL}) more explicitly.

\begin{prop}
The commutation relations can be written in the form
\begin{eqnarray}
B^{(1)}(z_1) \ B^{(2)}(z_2) & = & B^{(2)}(z_2) \ B^{(1)}(z_1) \ 
\Rred^{(12)}(z_1-z_2, \lared), \nonumber \\
A(z_1) \ B(z_2) & = & B(z_2) \ A(z_1) \ a(z_2-z_1,\lambda) + B(z_1) 
\ A(z_2) \ b(z_2-z_1,\lambda), \nonumber \\
\label{comm} 
D^{(1)}(z_1) \ B^{(2)}(z_2) & = & a^{(1)}(z_1-z_2, \lambda - \gamma h) \ 
B^{(2)}(z_2) \ D^{(1)}(z_1) \ \Rred^{(12)}(z_1-z_2, \lared) \\
& & \hspace*{20pt} + \ c^{(1)}(z_1-z_2, \lambda - \gamma h) \ B^{(2)}(z_1) 
\ D^{(1)}(z_2) \ \Pred^{(12)}, \nonumber 
\end{eqnarray} 
where $\lared=(\lambda_2, \ldots, \lambda_N)$, 
$\Rred$ is the $R$-matrix associated to $\mathrm{gl}_{N-1}$
\begin{displaymath}
\Rred(z,\lared) \ = \ \sumi{2} e_{ii} \otimes e_{ii} + \sumij{2} 
\alpha(z,\lambda_i - \lambda_j) \ e_{ii} \otimes e_{jj} 
+ \sumij{2} \beta(z,\lambda_i - \lambda_j) \ e_{ij} \otimes e_{ji},
\end{displaymath}
$\Pred^{(12)}=\sum_{i,j=2}^N e_{ij} \otimes e_{ji}$ is the 
permutation operator, and 
$a(z,\lambda)$, $b(z,\lambda)$, $c(z,\lambda)$ are diagonal 
$(N-1) \times (N-1)$ matrices with elements $(j=2, \ldots, N)$
\begin{displaymath}
[a(z,\lambda)]_{jj} = \frac{1}{\alpha(z,\lambda_j-\lambda_1)}, \ \ \ \  
[b(z,\lambda)]_{jj} = - \ \frac{\beta(z,\lambda_1-\lambda_j)}
{\alpha(z,\lambda_j-\lambda_1)}, \ \ \ \ 
[c(z,\lambda)]_{jj} = - \ \frac{\beta(z,\lambda_j-\lambda_1)}
{\alpha(z,\lambda_j-\lambda_1)}.
\end{displaymath}
\end{prop}

\begin{myproof}
If we denote by $R_{ik,jn}$ the elements of the $R$-matrix:
\begin{displaymath}
R^{(12)}(z,\lambda) \ = \ \sum_{i,j,k,n=1}^N \ R_{ik,jn}(z,\lambda) \ 
e_{ik} \otimes e_{jn}, 
\end{displaymath}
the element $(ik,jn)$ of the matrix relation~(\ref{RLL}) reads
\begin{displaymath}
\sum_{r,s=1}^N \ R_{ir,js}(z_-w,\lambda-\gamma h) \ 
\LL_{rk}(z) \ \LL_{sn}(w) \ = \ 
\sum_{r,s=1}^N \ \LL_{js}(w) \ \LL_{ir}(z) \ R_{rk,sn}(z-w,\lambda). 
\end{displaymath}
Taking successively $i=j=1$ and $k,n \geq 2$; $i=j=n=1$ and 
$k \geq 2$; $j=1$ and $i,k,n \geq 2$, and using the particular form 
of the $R$-matrix ($R_{ik,jn}=0$ unless $(k,n)=(i,j)$), we get
\begin{eqnarray*}
& & B_k(z) \ B_n(w) \ = \ \sum_{r,s=2}^N \ B_s(w) \ B_r(z) \ 
R_{rk,sn}(z-w,\lambda), \\ 
& & B_k(z) \ A(w) \ = \ B_k(w) \ A(z) \ \beta(z-w,\lambda_1-\lambda_k) \ 
+ \ A(w) \ B_k(z) \ \alpha(z-w,\lambda_k-\lambda_1), \\ 
& & \vphantom{\sum_{r,s=2}^N} 
\beta(z-w,(\lambda-\gamma h)_i-(\lambda-\gamma h)_1) \ B_k(z) \ 
D_{in}(w) \ + \ \alpha(z-w,(\lambda-\gamma h)_i-(\lambda-\gamma h)_1) \ 
D_{ik}(z) \ B_n(w) \\ 
& & \hspace*{20pt} = \ \sum_{r,s=2}^N \ B_s(w) \ D_{ir}(z) \ 
R_{rk,sn}(z-w,\lambda),
\end{eqnarray*}
which can be written in the form~(\ref{comm}).
\end{myproof}

Since the coefficients appearing in the commutation relations~(\ref{comm}) 
depend on $\lambda$ and $h$, we also need to know the way the functions 
$f(\lambda,h)$ go through the different operators. 

\begin{prop}
For any complex-valued function $f$ of $\lambda$ and $h$, we have:
\begin{eqnarray}
\nonumber
A(z) \ f(\lambda,h) & = & f(\lambda-\gamma\omega_1,h) \ A(z), \\
\label{ABDshift}
B_i(z) \ f(\lambda,h) & = & f(\lambda-\gamma\omega_i,h+\omega_1-\omega_i) 
\ B_i(z), \ \ \ \ \mathrm{for} \ i=2, \ldots, N, \\
\nonumber
D_{ij}(z) \ f(\lambda,h) & = & f(\lambda-\gamma\omega_j,h+\omega_i-\omega_j) 
\ D_{ij}(z), \ \ \ \ \mathrm{for} \ i,j=2, \ldots, N.
\end{eqnarray}
\end{prop}

\begin{myproof}
For the $\lambda$ dependence it is a straightforward consequence 
of the definition ~(\ref{defL}) of the Lax operator. For the $h$ 
dependence, it is not difficult to see that since $L$ commutes 
with the action of $\h$, we also have 
$[\LL^{(12)}(z) , f(h^{(1)}+h^{(2)}) ] = 0 $, or equivalently 
$\LL_{ij}(z) f(h+\omega_j) = f(h+\omega_i) \LL_{ij}(z)$.
\end{myproof}

\section{The first step of the Bethe ansatz}

\label{5}

With the help of the commutation relations~(\ref{comm}, \ref{ABDshift}), 
the action of $\T=A+\tr D$ on the vector $\psi$ given by~(\ref{psi}) 
can be recast in the form
\begin{eqnarray*}
& & \T(z) \ \psi(t_1,\ldots,t_m,\lambda) \ = \ B^{(1)}(t_1) \cdots B^{(m)}
(t_m) \ v \ g_0^{(1,\ldots,m)}(\lambda) \\
& & \hspace*{20pt} + \sum_{k=1}^m 
B^{(k)}(z) B^{(k+1)}(t_{k+1}) \cdots B^{(m)}(t_m) B^{(1)}(t_1) \cdots 
B^{(k-1)}(t_{k-1}) \ v \ g_k^{(1,\ldots,m)}(\lambda).
\end{eqnarray*}
Here and until the end of this section, we write $m \equiv (N-1)\ell$. 
Admitting that these different terms are linearly independent, we have 
to impose that $g_k=0$ for $k=1,\ldots,m$, and that $g_0$ is proportional 
to $g$. The term associated to $g_0$ is usually called the ``wanted term'' 
and  the other ones the ``unwanted terms''. 

If we carried out the commutation procedure, we would get $2^m$ terms 
for the $A$ part, and as many for the $D$ part: we see that this 
is hopeless to apply the process literally for all the terms. Nevertheless 
the ``wanted term'' and the first ``unwanted term'' are easy to obtain. 

To get the ``wanted term'', we  need to keep only the first 
part of the commutation relations $A(z) B(t) = B(t) A(z) a(t-z) + 
(\cdots)$, $D^{(1)}(z) B^{(2)}(t) = a^{(1)}(z-t) B^{(2)}(t) D^{(1)}(z) 
\Rred^{(12)}(z-t) + (\cdots)$. Between two commutations of $A(z)$ with one 
$B(t)$, we push a factor $a(t-z)$ to the right, and between two commutations 
of $D(z)$ with one $B(t)$, we push a factor $\Rred$ to the right, using the 
relations ~(\ref{ABDshift}). Once $A(z)$ and $D(z)$ are completely to 
the right, they act on $v \ g(\lambda)$ according to ~(\ref{Aaction}, 
\ref{Daction}). We get eventually 
\begin{eqnarray*}
\lefteqn{g_0^{(1,\ldots,m)}(\lambda) \ = \ 
\left[ \vphantom{\sum_{k=2}^m} a^{(m)}\Big(t_m-z,\lambda_1-\gamma,
\lared\Big) \ a^{(m-1)}\Big(t_{m-1}-z,\lambda_1-\gamma,\lared+\gamma 
\hred^{(m)}\Big) \ \cdots \right.} \\
& & \cdots \ \left. a^{(1)}\Big(t_1-z,\lambda_1-\gamma,\lared+\gamma 
\sum_{k=2}^m \hred^{(j)}\Big) \right] \ 
g^{(1,\ldots,m)}(\lambda_1-\gamma,\lared) \\
& & + \ \tr^{(0)} \left[ \vphantom{\sum_{k=2}^m} 
a^{(0)}(z-t_1,\lambda_1-\ell\gamma,\lared) \  
a^{(0)}(z-t_2,\lambda_1-(\ell+1)\gamma,\lared) \ \cdots \right. \\
& & \vphantom{\sum_{k=2}^m} \cdots \ 
a^{(0)}(z-t_m,\lambda_1-(\ell+m-1)\gamma,\lared) \ 
\Phi^{(0)}(z) \ \Rred^{(0m)}\Big(z-t_m,\lared\Big) \\ 
& & \left. 
\Rred^{(0,m-1)}\Big(z-t_{m-1},\lared+\gamma \hred^{(m)}\Big) \ \cdots \ 
\Rred^{(0,1)}\Big(z-t_1,\lared+\gamma \sum_{k=2}^m\hred^{(k)}\Big) 
\right] g^{(1,\ldots,m)}(\lambda).
\end{eqnarray*}
What we have in mind is to be sent back to the diagonalization 
of a transfer matrix $\Tred$ of the same kind as $\T$, with dimension $N$ 
decreased by 1. This is the case if we choose the right dependence of the 
vector $g$ in the variable $\lambda_1$.

\begin{prop}
If we take $g^{(1,\ldots,m)}(\lambda) = G(\lambda) \ \psired(\lared)$ with 
\begin{equation}
\label{G}
G(\lambda) \ = \ \mathrm{e}^{c_1 \lambda_1} \prod_{j=2}^N \prod_{p=1}^\ell 
\theta(\lambda_1-\lambda_j-p\gamma),
\end{equation}
where $c_1$ is an arbitrary constant, then
\begin{eqnarray}
\label{wanted}
\lefteqn{g_0^{(1,\ldots,m)}(\lambda) \ = \ 
G(\lambda) \left[ \vphantom{\frac{\theta(z)}{\theta(z-N\ell\gamma)}} 
S_+(z;t_1,\ldots,t_m) \ \mathrm{e}^{-\gamma c_1} \right.} \\
& & \left. + \ \frac{\theta(z)}{\theta(z-N\ell\gamma)} \ 
S_-(z;t_1,\ldots,t_m) \ \Tred^{(m,\ldots,1)}(z) \right] \ \psired(\lared), 
\nonumber
\end{eqnarray}
where $ \displaystyle 
S_\pm(z;t_1,\ldots,t_m) = \prod_{k=1}^{m} \frac{\theta(z-t_k\pm\gamma)}
{\theta(z-t_k)}$ and
\begin{eqnarray*}
\lefteqn{\Gammared^{(1,\ldots,m)} \ \Tred^{(m,\ldots,1)}(z) \ 
\left(\Gammared^{(1,\ldots,m)}\right)^{-1}} \\
& & = \tr^{(0)} \left[
\Rred^{(0m)}\Big(z-t_m,\lared-\gamma\sum_{j=1}^{m-1} \hred^{(j)}\Big) 
\ \cdots \ \Rred^{(01)}\Big(z-t_1,\lared\Big) \ \Gammared^{(0)} \right].
\end{eqnarray*}
Up to a conjugation by 
$\Gammared^{(1,\ldots,m)} = (\Gammared^{(1)} \ \cdots \ 
\Gammared^{(m)})$ (where $\Gammared$ is the shift operator 
$\Gammared=\mathrm{diag}(\Gamma_i)_{i=2,\ldots,N}$), 
$\Tred$ is the transfer matrix acting in 
$\Wred=(\C^{N-1})^{\otimes (N-1)\ell}$ given by an expression 
similar to ~(\ref{genT}, \ref{genL}).
\end{prop}

\begin{myproof} 
Let us look at the first term of $g_0$, which is given by 
the action of a diagonal matrix on $g$. Because $g$ is a zero weight 
vector, we find that this term is simply a multiple of $g$ given by
\begin{eqnarray*}
\lefteqn{\left[ a^{(m)} \ a^{(m-1)} \ \cdots \ a^{(1)} 
\right] g^{(1,\ldots,m)}(\lambda-\gamma\omega_1)} \\ 
& & = \ \prod_{k=1}^m 
\frac{\theta(z-t_k+\gamma)}{\theta(z-t_k)} \ \prod_{j=2}^N 
\frac{\theta(\lambda_1-\lambda_j-\gamma)}{\theta(\lambda_1-\lambda_j-
(\ell+1)\gamma)} \ g^{(1,\ldots,m)}(\lambda-\gamma\omega_1).
\end{eqnarray*}
For this to be exactly $g^{(1,\ldots,m)}(\lambda)$ up to a constant 
factor, we have to take 
\begin{displaymath}
g^{(1,\ldots,m)}(\lambda) \ = \ G(\lambda) \ \psired(\lared)
\end{displaymath}
with the function $G$ as in ~(\ref{G}). 

Now let us simplify the second term. If we compute the $j$th element of 
the matrix $a(z-t_1,\lambda_1-\ell\gamma,\lared) \cdots 
a(z-t_m,\lambda_1-(\ell+m-1)\gamma,\lared) \Phi(z)$, we find 
\begin{eqnarray*}
\lefteqn{\left[ a(z-t_1,\lambda_1-\ell\gamma,\lared) \ \cdots \ 
a(z-t_m,\lambda_1-(\ell+m-1)\gamma,\lared) \ \Phi(z)\right]_j} \\
& & = \ \frac{\theta(z)}{\theta(z-(m+\ell)\gamma)} \ \prod_{k=1}^m 
\frac{\theta(z-t_k-\gamma)}{\theta(z-t_k)} \ 
\frac{\theta(\lambda_1-\lambda_j-\ell\gamma)}{\theta(\lambda_1-\lambda_j)} \ 
\Gamma_j
\end{eqnarray*}
and then it is clear that  
\begin{eqnarray*}
\lefteqn{\left[ a(z-t_1,\lambda_1-\ell\gamma,\lared) \ \cdots \ 
a(z-t_m,\lambda_1-(\ell+m-1)\gamma,\lared) \ \Phi(z)\right]_j \ G(\lambda)} \\
& & = \ \frac{\theta(z)}{\theta(z-(m+\ell)\gamma)} \ \prod_{k=1}^m 
\frac{\theta(z-t_k-\gamma)}{\theta(z-t_k)} \ G(\lambda) \  \Gamma_j.
\end{eqnarray*}
So when we take the trace we simply get the transfer 
matrix $\Tred$ up to a scalar coefficient.
\end{myproof}

The first ``unwanted term'' is obtained similarly, except that 
for the commutation of $A(z)$ and $D(z)$ accross $B(t_1)$, we keep the 
second part of the commutation relations:
\begin{eqnarray*}
\lefteqn{g_1^{(1,\ldots,m)}(\lambda) \ = \ 
\left[ \vphantom{\sum_{k=2}^m} a^{(m)}\Big(t_m-t_1,\lambda_1-\gamma,
\lared\Big) \ \cdots \ a^{(2)}\Big(t_2-t_1,\lambda_1-\gamma,
\lared+\gamma \sum_{k=3}^m \hred^{(k)}\Big) \right.} \\
& & \ \left. 
\ b^{(1)}\Big(t_1-z,\lambda_1-\gamma,\lared+\gamma 
\sum_{k=2}^m \hred^{(k)}\Big) \right] \ 
g^{(1,\ldots,m)}(\lambda_1-\gamma,\lared) \\
& & + \ \tr^{(0)} \left[ \vphantom{\sum_{k=2}^m} 
c^{(0)}(z-t_1,\lambda_1-\ell\gamma,\lared) \ 
a^{(0)}(t_1-t_2,\lambda_1-(\ell+1)\gamma,\lared) \ \cdots \right. \\ 
& & \vphantom{\sum_{k=2}^m} \cdots \ 
a^{(0)}(t_1-t_m,\lambda_1-(\ell+m-1)\gamma,\lared) \ 
\Phi^{(0)}(t_1) \ \Rred^{(0m)}\Big(t_1-t_m,\lared\Big) \ \cdots \\
& & \left. \cdots \ 
\Rred^{(0,2)}\Big(t_1-t_2,\lared+\gamma \sum_{k=3}^m\hred^{(k)}\Big) \ 
\Pred^{(01)} \right] g^{(1,\ldots,m)}(\lambda).
\end{eqnarray*}

\begin{prop} With $g^{(1,\ldots,m)}(\lambda) = G(\lambda) \ 
\psired(\lared)$ where $G(\lambda)$ is given 
by~(\ref{G}), we have 
\begin{eqnarray*}
\lefteqn{g_1^{(1,\ldots,m)}(\lambda) \ = \ G(\lambda) \ 
\frac{\theta(\gamma)}{\theta(t_1-z)} \ X^{(1)}(z,t_1,\lambda) \ 
\left[  \vphantom{\frac{\theta(t_1)}{\theta(t_1-N\ell\gamma)}} \ 
K_+(t_1,\ldots,t_m) \ \mathrm{e}^{-\gamma c_1} \right.} \\
& & \left. - \ \frac{\theta(t_1)}{\theta(t_1-N\ell\gamma)} \ 
K_-(t_1,\ldots,t_m) \ \Tred^{(m,\ldots,1)}(t_1) \right] \ \psired(\lared)
\end{eqnarray*}
where
\begin{displaymath}
K_\pm(t_1,\ldots,t_m) = \prod_{k=2}^{m} 
\frac{\theta(t_1-t_k\pm\gamma)}{\theta(t_1-t_k)}
\end{displaymath}
and $X$ is the diagonal $(N-1) \times (N-1)$ matrix of elements 
$ \displaystyle 
[X(z,t,\lambda)]_{jj} = \frac{\theta(z-t+\lambda_j-\lambda_1+\ell\gamma)}
{\theta(\lambda_j-\lambda_1+\ell\gamma)}$, for $j=2,\ldots,N$. 
\end{prop}

We shall not compute directly the other ``unwanted terms''. 
To obtain them, let us change the order of 
the $B$ operators in ~(\ref{psi}), using the commutation 
relations~(\ref{comm},\ref{ABDshift}). 
We get 
\begin{equation}
\label{newpsi}
\psi(t_1,\ldots,t_m,\lambda) = B^{(2)}(t_2) \cdots B^{(m)}(t_m) 
B^{(1)}(t_1)\ v \ \Rcalred^{(1;m,\ldots,2)} \ g^{(1,\ldots,m)}(\lambda), 
\end{equation}
where
\begin{displaymath}
\Rcalred^{(1;m,\ldots,2)} \ = \ 
\Rred^{(1m)}\Big(t_1-t_m,\lared\Big) \ \cdots \ 
\Rred^{(12)}\Big(t_1-t_2,\lared+\gamma\sum_{k=3}^m \hred^{(k)}\Big).
\end{displaymath}

Starting with~(\ref{newpsi}) and 
applying the same procedure as before, we find another expression 
for the ``wanted term''; let us check that the two results are equivalent. 
We find now
\begin{eqnarray*}
\lefteqn{g_0^{(1,\ldots,m)}(\lambda) \ = \ G(\lambda) \ 
\left( \Rcalred^{(1;m,\ldots,2)} \right)^{-1}
\left[ \vphantom{\frac{\theta(z)}{\theta(z-N\ell\gamma)}} 
S_+(z;t_2,\ldots,t_m,t_1) \ \mathrm{e}^{-\gamma c_1} \right.} \\
& & \left. + \ \frac{\theta(z)}{\theta(z-N\ell\gamma)} \ 
S_-(z;t_2,\ldots,t_m,t_1) \ \Tred^{(1,m,\ldots,2)}(z) 
\right] \ \Rcalred^{(1;m,\ldots,2)} \ \psired(\lared).
\end{eqnarray*}
Since $S_\pm$ are symmetric functions of $(t_1,\ldots,t_m)$, and thanks 
to some commutation property of $\Tred$ and $\Rcalred$ which is written 
below, it is straightforward to see that this second expression is 
equal to~(\ref{wanted}). 

\begin{prop} 
\begin{equation}
\label{TR}
\Tred^{(1,m,\ldots,2)}(z) \ \Rcalred^{(1;m,\ldots,2)} \ = 
\ \Rcalred^{(1;m,\ldots,2)} \ \Tred^{(m,\ldots,1)}(z).
\end{equation}
\end{prop}

\begin{myproof} 
We can write 
$\Tred^{(m,\ldots,1)}= \tr^{(0)} 
\tilde{\LL}^{(0;m,\ldots,1)}$ with 
\begin{eqnarray*}
\tilde{\LL}^{(0;m,\ldots,1)} & = & \Gaminv{m} \RR{0}{m} \Gaminv{m-1} 
\RR{0}{m-1} \cdots \Gaminv{1} \RR{0}{1} \Gam{0} \Gam{1} \cdots \Gam{m} \\
& = & \Gam{0} \RR{0}{m} \Gaminv{m} \RR{0}{m-1} \Gaminv{m-1} \cdots 
\RR{0}{1} \Gam{2} \cdots \Gam{m}
\end{eqnarray*}
because $\Gam{0} \Gam{j}$ commutes with $\RR{0}{j}$. 
Using the dynamical Yang--Baxter equation~(\ref{DYB}) in the form
\begin{displaymath}
\Gam{0} \RR{0}{1} \Gaminv{1} \RR{0}{j} \Gaminv{j} \RR{1}{j} = \RR{1}{j} 
\Gaminv{j} \RR{0}{j} \Gam{0} \RR{0}{1} \Gaminv{1},
\end{displaymath}
one can see that $\tilde{\LL}^{(0;1,2)} \Rcalred^{(1;2)} = 
\Rcalred^{(1;2)} \tilde{\LL}^{(0;2,1)}$. Then using the relations 
\begin{eqnarray*}
\Rcalred^{(1;m+1,\ldots,2)} & = & \RR{1}{m+1} \Gaminv{m+1} 
\Rcalred^{(1;m,\ldots,2)} \Gam{m+1} \\
& = & \Gaminv{1} \Gaminv{m+1} \RR{1}{m+1} \Gam{1} 
\Rcalred^{(1;m,\ldots,2)} \Gam{m+1},
\end{eqnarray*}
it is easy to prove the result by recursion on $m$.
\end{myproof}

Let us note that at this point in the usual nested Bethe ansatz, the 
situation is a little simpler. Indeed if there is no shift entering 
the definition of the Lax operator, by cyclicity of the trace 
$\Tred^{(1,m,\ldots,2)}(z)=\Tred^{(m,\ldots,1)}(z)$, and 
$\Rcalred^{(1;m,\ldots,2)}$ is precisely equal to 
$\Tred^{(1,m,\ldots,2)}(t_1)$. In that case, the relation~(\ref{TR}) is 
just the expression of the commutation of the transfer matrices for 
different values of the spectral parameter. But here because of the 
dynamical feature of the Lax operator, $\Rcalred^{(1;m,\ldots,2)}$ is 
not a transfer matrix, and $\T^{(1,m,\ldots,2)}(z)$ is not equal to 
$\T^{(m,\ldots,1)}(z)$.

It is now clear that the expression of the second ``unwanted term'' is:
\begin{eqnarray*}
\lefteqn{g_2^{(1,\ldots,m)}(\lambda) \ = \ 
G(\lambda) \ \frac{\theta(\gamma)}{\theta(t_2-z)} \ X^{(2)}(z,t_2,\lambda) 
\left[ \vphantom{\frac{\theta(t_2)}{\theta(t_2-N\ell\gamma)}} 
K_+(t_2,\ldots,t_m,t_1) \ \mathrm{e}^{-\gamma c_1} \right.} \\
& & \left. - \ \frac{\theta(t_2)}{\theta(t_2-N\ell\gamma)} \ 
K_-(t_2,\ldots,t_m,t_1) \ 
\Tred^{(1,m,\ldots,2)}(t_2) \right] 
\Rcalred^{(1;m,\ldots,2)} \ \psired(\lared),
\end{eqnarray*}
which, with the use of relation~(\ref{TR}),  can be written
\begin{eqnarray*}
\lefteqn{g_2^{(1,\ldots,m)}(\lambda) \ = \ 
G(\lambda) \ \frac{\theta(\gamma)}{\theta(t_2-z)} \ X^{(2)}(z,t_2,\lambda) 
\ \Rcalred^{(1;m,\ldots,2)}} \\
& & \left[ K_+(t_2,\ldots,t_m,t_1) \ \mathrm{e}^{-\gamma c_1} \ 
- \ \frac{\theta(t_2)}{\theta(t_2-N\ell\gamma)} \ K_-(t_2,\ldots,t_m,t_1) 
\ \Tred^{(m,\ldots,1)}(t_2)  \right] \ \psired(\lared).
\end{eqnarray*}

The expressions of the other ``unwanted terms'' are obtained by 
repeating several times the permutation of the $B$'s. 

\begin{prop} The cancellation of all the ``unwanted terms'' is 
equivalent to the set of relations  
\begin{displaymath}
\Tred^{(m,\ldots,1)} (t_k) \ \psired(\lared) \ = \ 
\mathrm{e}^{-\gamma c_1} \ 
\frac{\theta(t_k-N\ell\gamma)}{\theta(t_k)} \  \prod_{^{i=1}_{i\neq k}}^{m}
\frac{\theta(t_k-t_i+\gamma)}{\theta(t_k-t_i-\gamma)}
\ \psired(\lared), \ \ \ \ \forall k=1,\ldots, m.
\end{displaymath}
\end{prop}

\section{The Bethe equations}

\label{6}

After the first step of the Bethe ansatz described in the preceding 
section, we are led to the problem of the diagonalization of 
$\Tred^{((N-1)\ell,\ldots,1)}(z)$ 
for arbitrary values of the parameters $t_i$. 
Let us call $\psired^{((N-1)\ell,\ldots,1)}(\lared)$ 
the zero weight eigenstate of 
$\Tred^{((N-1)\ell,\ldots,1)}(z)$, 
with eigenvalue $\epsred(z)$ (both $\epsred$ and $\psired$ 
depending implicitely on the parameters $t_i$). If the 
following conditions are satisfied:
\begin{displaymath}
\epsred (t_k) \ = \ \mathrm{e}^{-\gamma c_1} \ 
\frac{\theta(t_k-N\ell\gamma)}{\theta(t_k)} \  \prod_{^{i=1}_{i\neq k}}
^{(N-1)\ell} 
\frac{\theta(t_k-t_i+\gamma)}{\theta(t_k-t_i-\gamma)}, \ \ \ \ 
k=1,\ldots, (N-1)\ell, 
\end{displaymath}
the vector $\psi(t_1,\ldots,t_{(N-1)\ell},\lambda)$ given by~(\ref{psi}) 
will be an eigenstate of $\T(z)$ with eigenvalue
\begin{displaymath}
\varepsilon(z) \ = \ \prod_{k=1}^{(N-1)\ell} \frac{\theta(z-t_k+\gamma)}
{\theta(z-t_k)} \ \mathrm{e}^{-\gamma c_1} \ + \ 
\frac{\theta(z)}{\theta(z-N\ell\gamma)} \ 
\prod_{k=1}^{(N-1)\ell} \frac{\theta(z-t_k-\gamma)}{\theta(z-t_k)} \ 
\epsred(z).
\end{displaymath}

The second step of the Bethe ansatz thus consists in diagonalizing 
$\Tred^{((N-1)\ell,\ldots,1)}(z)$ by repeating 
the same procedure. 
The total shift $\Gammared^{(1,\ldots,(N-1)\ell)}$ entering the 
definition of $\Tred^{((N-1)\ell,\ldots,1)}$ has no action because 
we are looking for zero weight eigenstates only. 
We seek eigenstates in the form 
\begin{displaymath}
\psired(u_1,\ldots,u_{(N-2)\ell},\lared) = 
\Bred^{(1)}(u_1) \cdots \Bred^{((N-2)\ell)}
(u_{(N-2)\ell}) \ \vred \ \gred^{(1,\ldots,(N-2)\ell)}(\lared), 
\end{displaymath}
where $\vred$ is the vector of $\Wred$ 
given by $\vred = e_2 \otimes \cdots \otimes e_2$ 
and $\Bred_j=\Lred_{2j}$ for $j=3,\ldots,N$. Everything works out 
similarly, except the fact that the operators 
$\Dred_{ij}=\Lred_{ij}$ for $i,j=3,\ldots,N$ act a little differently 
on $\vred$ because of the parameters $t_k$ in $\Tred$:
\begin{displaymath}
\Dred(z) \ [g(\lared) \ \vred] \ = \ 
\sumi{3} e_{ii} \ \Phired_i(z,\lared) \ g(\lared-\omegared_i) \ \vred, 
\end{displaymath}
where $ \displaystyle 
\Phired_i(z,\lared) = \frac{\theta(\lambda_i-\lambda_2+(N-1)\ell\gamma)}
{\theta(\lambda_i-\lambda_2)} \ \prod_{k=1}^{(N-1)\ell} 
\frac{\theta(z-t_k)}{\theta(z-t_k-\gamma)}$.

We continue in this way until the $N$th step; at this last step, 
the diagonalization problem to solve is simply 
$\psi(\lambda_N-\gamma) = \varepsilon_N \ \psi(\lambda_N)$, 
the solution of which is (up to a constant) 
$\psi(\lambda_N) = \mathrm{e}^{c_N \lambda_N}, 
\varepsilon_N = \mathrm{e}^{-\gamma c_N}$. 

The whole procedure can be summarized as follows.

\begin{mydef} 
We introduce a set of parameters 
$(t^{(n)}_1,\ldots,t^{(n)}_{(N-n)\ell})$ for $n=0,\ldots,N-1$, with
\begin{displaymath}
t^{(0)}_j=(j-1)\gamma, \ \ \ \ \forall j=1,\ldots,N\ell.
\end{displaymath}
We define a set of functions $\varepsilon_n(z)$, $n=1,\ldots,N$, 
by the recursion relation 
\begin{eqnarray*}
\varepsilon_N(z) & = & \mathrm{e}^{-\gamma c_N}, \\
\varepsilon_n(z) & = & \mathrm{e}^{-\gamma c_n} \prod_{k=1}^{(N-n)\ell} 
\frac{\theta(z-t_k^{(n)}+\gamma)}{\theta(z-t_k^{(n)})} 
\ + \prod_{i=1}^{(N-n+1)\ell} 
\frac{\theta(z-t_i^{(n-1)})}{\theta(z-t_i^{(n-1)}-\gamma)} 
\prod_{k=1}^{(N-n)\ell} 
\frac{\theta(z-t_k^{(n)}-\gamma)}{\theta(z-t_k^{(n)})} 
\varepsilon_{n+1}(z), 
\end{eqnarray*}
where $c_1,\ldots,c_N$ are arbitrary complex parameters. 

For each $n=1,\ldots,N-1$, we define a Lax matrix by setting:
\begin{eqnarray*}
& & R\left(n\big|z,\lambda\right) \ = \ \sumi{n} e_{ii} \otimes e_{ii} 
+ \sumij{n} \alpha(z,\lambda_i-\lambda_j) \ e_{ii} \otimes e_{jj} 
+ \sumij{n} \beta(z,\lambda_i-\lambda_j) \ e_{ij} \otimes e_{ji}, \\
& & \LL\left(n\big|z\right) \ = \ \left(
\prod_{k=1}^{(N-n+1)\ell} 
R^{(0k)}\Big(n\big|z-t_k^{(n-1)},\lambda-\gamma
\sum_{j \bowtie k} h^{(j)} \Big) \right) \ \Gamma(n)^{(0)},
\end{eqnarray*}
where the product on $k$ is taken in increasing order if $n$ is odd, and 
decreasing order if $n$ is even. The notation $j \bowtie k$ means 
$j>k$ if $n$ is odd, and $j<k$ if $n$ is even. $\Gamma(n)$ is the 
matrix $\mathrm{diag}(\Gamma_i)_{i=n,\ldots,N}$. 
The creation operators are the elements of $\LL$ given by 
\begin{displaymath}
B_j\left(n\big|t\right) \ = \ \LL_{nj}\left(n\big|t\right), \ \ \ \ 
j=n+1,\ldots,N,
\end{displaymath}
where the subscripts $nj$ refer to the auxiliary space (0).

Finally we define the states $\psi_n(\lambda_n,\ldots,\lambda_N)$ by 
$\displaystyle \psi_N(\lambda_N) = \mathrm{e}^{c_N \lambda_N}$ 
and the recursion relation (for $n=1,\ldots,N-1$)
\begin{displaymath}
\psi_n \ = \ 
B^{(1)}\left(n\big|t_1^{(n)}\right) \ \cdots \ 
B^{((N-n)\ell)}\left(n\big|t_{(N-n)\ell}^{(n)}\right) \ 
v_n \ G_n \ \psi_{n+1},
\end{displaymath}
where 
\begin{eqnarray*}
& & v_n \ = \ e_n\otimes \cdots \otimes e_n \ \in 
(\C^{N-n+1})^{\otimes (N-n+1)\ell}, \\
& & G_n(\lambda_n,\ldots,\lambda_N) \ = \ \mathrm{e}^{c_n \lambda_n} \ 
\prod_{j=n+1}^{N} \prod_{p=1}^{\ell} \ \theta(\lambda_n-\lambda_j-p\gamma). 
\end{eqnarray*}
\end{mydef}

\begin{prop}
If the parameters $t_k^{(n)}$ are solution of the Bethe equations 
\begin{eqnarray}
& & \prod_{^{j=1}_{j\neq k}}^{(N-n)\ell} 
\frac{\theta(t_k^{(n)}-t_j^{(n)}+\gamma)}
{\theta(t_k^{(n)}-t_j^{(n)}-\gamma)} \ 
\prod_{i=1}^{(N-n+1)\ell} 
\frac{\theta(t_k^{(n)}-t_i^{(n-1)}-\gamma)}{\theta(t_k^{(n)}-t_i^{(n-1)})} \ 
\prod_{i=1}^{(N-n-1)\ell} 
\frac{\theta(t_k^{(n)}-t_i^{(n+1)})}{\theta(t_k^{(n)}-t_i^{(n+1)}+\gamma)} \\
& & \hspace*{20pt} = \ \mathrm{e}^{\gamma(c_n-c_{n+1})}, \ \ \ \ \ \ \ \ 
\forall n=1,\ldots,N-1, \ \forall k=1,\ldots,(N-n)\ell,
\nonumber
\end{eqnarray} 
then $\psi_1(\lambda_1,\ldots,\lambda_N)$ is a 
zero weight eigenstate of $\T(z)=\tr^{(0)} \LL\left(1\big|z\right)$ 
with eigenvalue 
\begin{displaymath}
\varepsilon_1(z) = \sum_{i=1}^N \mathrm{e}^{-\gamma c_i} \ 
\prod_{m=1}^{(N-i)\ell} 
\frac{\theta(z-t_m^{(i)}+\gamma)}{\theta(z-t_m^{(i)})} \ 
\prod_{j=1}^{i-1} \ \left( \prod_{m=1}^{(N-j)\ell} 
\frac{\theta(z-t_m^{(j)}-\gamma)}{\theta(z-t_m^{(j)})} \ 
\prod_{m=1}^{(N-j+1)\ell} 
\frac{\theta(z-t_m^{(j-1)})}{\theta(z-t_m^{(j-1)}-\gamma)} \right).
\end{displaymath} 
\end{prop}

\begin{myproof}
The cancellation of all the ``unwanted terms'' is equivalent to the 
set of equations 
\begin{displaymath}
\varepsilon_{n+1}(t_k^{(n)}) \ = \ \mathrm{e}^{-\gamma c_n} \ 
\prod_{i=1}^{(N-n+1)\ell} 
\frac{\theta(t_k^{(n)}-t_i^{(n-1)}-\gamma)}{\theta(t_k^{(n)}-t_i^{(n-1)})} 
\ \prod_{^{j=1}_{j\neq k}}^{(N-n)\ell} 
\frac{\theta(t_k^{(n)}-t_j^{(n)}+\gamma)}
{\theta(t_k^{(n)}-t_j^{(n)}-\gamma)}, 
\end{displaymath} 
for $n=1,\ldots,N-1$ and $ k=1,\ldots,(N-n)\ell$. It is easy to see that 
\begin{eqnarray*}
\varepsilon_n(z) = \sum_{i=n}^N \mathrm{e}^{-\gamma c_i} \ 
\prod_{m=1}^{(N-i)\ell} 
\frac{\theta(z-t_m^{(i)}+\gamma)}{\theta(z-t_m^{(i)})} \ 
\prod_{j=n}^{i-1} \ \left( \prod_{m=1}^{(N-j)\ell} 
\frac{\theta(z-t_m^{(j)}-\gamma)}{\theta(z-t_m^{(j)})} \ 
\prod_{m=1}^{(N-j+1)\ell} 
\frac{\theta(z-t_m^{(j-1)})}{\theta(z-t_m^{(j-1)}-\gamma)} \right), 
\end{eqnarray*} 
and so in the expression of $\varepsilon_{n+1}(t_k^{(n)})$ there 
is only one non zero term, corresponding to $i=n+1$. 
\end{myproof}

\begin{remark}
Since the Ruijsenaars operator $M$ is related to the transfer matrix $\T(z)$ 
by the relation~(\ref{Ruij}), $\psi_1(\lambda_1,\ldots,\lambda_N)$ is 
an eigenfunction of $M$ with eigenvalue 
$\varepsilon = \frac{\theta(z-\gamma N\ell)}{\theta(z-\gamma\ell)} 
\varepsilon_1(z)$. This quantity does not depend on $z$, so 
we can evaluate it at $z=0$. We find
\begin{displaymath}
\varepsilon \ = \ \mathrm{e}^{-\gamma c_1} \ 
\frac{\theta(\gamma N\ell)}{\theta(\gamma\ell)} \ 
\prod_{m=1}^{(N-1)\ell} 
\frac{\theta(t_m^{(1)}-\gamma)}{\theta(t_m^{(1)})}.
\end{displaymath}
\end{remark}

\bigskip
\begin{center} {\sc Acknowledgement} \end{center}
I thank Giovanni Felder for discussions.


\bigskip 

\begin{thebibliography}{99}

\bibitem{FeVa97} G. Felder, A. Varchenko, {\it Elliptic quantum groups and 
Ruijsenaars models}, {\tt q-alg/9704005}.

\bibitem{FeVa96} G. Felder, A. Varchenko, {\it Algebraic Bethe ansatz 
for the elliptic quantum group $E_{\tau,\eta}(sl_2)$}, 
\np{B480}{1996}{485--503} and {\tt q-alg/9605024}. 

\bibitem{SkTaFa79} E.K. Sklyanin, L.A. Takhtadjan, L.D. Faddeev, {\it 
The quantum inverse problem method I}, \tmp{40}{1980}{688--706}.

\bibitem{Fa82} L.D. Faddeev, {\it Integrable models in $1+1$ dimensional 
quantum field theory}, Les Houches lectures 1982, \esp{1984}.

\bibitem{KuRe83} P.P. Kulish, N. Yu Reshetikhin, {\it Diagonalization of 
$gl(n)$ invariant transfer matrices and quantum $n$ wave system 
(Lee model)}, \jpamg{16}{1983}{L591--L596}.

\bibitem{FeVa95a} G. Felder, A. Varchenko, {\it Integral representations of 
solutions of the elliptic Knizhnik--Zamolodchikov--Bernard equations}, 
\imrn{5}{1995}{221--233} and {\tt hep-th/9502165}.

\bibitem{EtKi93a} P.I. Etingof, A.A. Kirillov, Jr., {\it Representations of 
affine Lie algebras, parabolic differential equations and Lam\'e 
functions}, \dmj{74}{1994}{585--614} and {\tt hep-th/9310083}.

\bibitem{EtKi94} P.I. Etingof, A.A. Kirillov, Jr., {\it On the affine 
analogue of Jack's and Macdonald's polynomials}, {\tt hep-th/9403168}.

\bibitem{ReVa94} N. Yu Reshetikhin, A. Varchenko, {\it Quasiclassical 
asymptotics of solutions to the KZ equations}, {\tt hep-th/9402126}.

\bibitem{FeVa95b} G. Felder, A. Varchenko, {\it Three formulae for 
eigenfunctions of integrable Schr\"odinger operators}, 
\cp{107}{1997}{143--175} and {\tt hep-th/9511120}.

\bibitem{EtKi93b} P.I. Etingof, A.A. Kirillov, Jr., {\it A unified 
representation theoretic approach to special functions}, 
\faa{28:1}{1994}{91--94} and {\tt hep-th/9312101}.

\bibitem{FetaVa96} G. Felder, V. Tarasov, A. Varchenko, {\it Solutions of 
the elliptic QKZB equations and Bethe ansatz I}, {\tt q-alg/9606005}.

\end{thebibliography}
\end{document}